\documentclass[11pt,twoside]{amsart}

\author{Valery Alexeev and Michel Brion}  
\address{Department of Mathematics\\
University of Georgia\\
Athens, GA 30602, USA}
\email{valery@math.uga.edu}
\address{Institut Fourier, B. P. 74\\
38402 Saint--Martin d'H\`eres Cedex, France}
\email{Michel.Brion@ujf-grenoble.fr}
\title{Stable reductive varieties II: Projective case}
\date{July 29, 2002; corrected version: April 15, 2003}

\theoremstyle{plain}
\newtheorem{theorem}{Theorem}[section]

\newcommand{\oM}{\overline{\operatorname{M}}}

\newcommand{\wgxg}{\wG\times\wG}
\newcommand{\uh}{\underline{h}}
\newcommand{\wY}{\widetilde{Y}}
\newcommand{\wlambda}{\tilde{\lambda}}
\newcommand{\wLambda}{\widetilde{\Lambda}}
\newcommand{\Chow}{\operatorname{Chow}}
\newcommand{\SP}{\operatorname{SP}}
\newcommand{\M}{\operatorname{M}}
\newcommand{\Mci}{{\operatorname{M}}^{\rm ci}}
\newcommand{\cMci}{{\mathcal M}^{\rm ci}}

\newcommand{\val}{\operatorname{val}}

\newcommand{\gxg}{G\times G}

\newcommand{\diag}{\operatorname{diag}}

\newcommand{\End}{\operatorname{End}}

\newcommand{\Conv}{\operatorname{Conv}}
\newcommand{\Aaut}{\operatorname{\bA ut}}
\newcommand{\gr}{\operatorname{gr}}

\newcommand{\bA}{{\mathbb A}}
\newcommand{\bC}{{\mathbb C}}

\newcommand{\bG}{{\mathbb G}}
\newcommand{\bK}{{\mathbb K}}
\newcommand{\bL}{{\mathbb L}}
\newcommand{\bM}{{\mathbb M}}
\newcommand{\bN}{{\mathbb N}}
\newcommand{\bP}{{\mathbb P}}
\newcommand{\bQ}{{\mathbb Q}}
\newcommand{\bR}{{\mathbb R}}

\newcommand{\bZ}{{\mathbb Z}}

\newcommand{\cD}{{\mathcal D}}

\newcommand{\cK}{{\mathcal K}}

\newcommand{\cM}{{\mathcal M}}

\newcommand{\cO}{{\mathcal O}}

\newcommand{\cR}{{\mathcal R}}

\newcommand{\hFun}{\widehat {\operatorname{Fun}}}

\newcommand{\inv}{^{-1}}

\newcommand{\hbK}{\widehat \bK}

\newcommand{\hbM}{\widehat \bM}

\newcommand{\wE}{\widetilde E}
\newcommand{\wG}{\widetilde G}

\newcommand{\wT}{\widetilde T}

\newcommand{\wX}{\widetilde X}
\newcommand{\wZ}{\widetilde Z}

\newcommand{\wDelta}{\widetilde\Delta}

\DeclareMathSymbol{\curvearrowright}{\mathrel}{AMSb}{"79}
\DeclareMathSymbol\rightsquigarrow {\mathrel}{AMSa}{"20}
\DeclareMathSymbol\square {\mathord}{AMSa}{"03}
\DeclareMathSymbol{\ltimes}         {\mathbin}{AMSb}{"6E}
\DeclareMathSymbol{\nmid}           {\mathrel}{AMSb}{"2D}
\DeclareMathSymbol{\twoheadrightarrow}  {\mathrel}{AMSa}{"10}

\newcommand{\ratmap}{- \kern -3pt \to}

\newcommand{\onto}{\twoheadrightarrow}

\newcommand{\follows}{\Rightarrow}

\newcommand{\lin}{\operatorname{lin}}

\newcommand{\Cone}{\operatorname{Cone}}
\newcommand{\Fun}{\operatorname{Fun}}

\newcommand{\id}{\operatorname{id}}

\newcommand{\Supp}{\operatorname{Supp}}

\newcommand{\Proj}{\operatorname{Proj}}

\newcommand{\Stab}{\operatorname{Stab}}
\newcommand{\Spec}{\operatorname{Spec}}

\newcommand{\Hom}{\operatorname{Hom}}

\newcommand{\Aut}{\operatorname{Aut}}

\newcommand{\chr}{\operatorname{char}}
\newcommand{\Sym}{\operatorname{Sym}}
\newcommand{\Ext}{\operatorname{Ext}}

\theoremstyle{plain}

\newtheorem{lemma}[theorem]{Lemma}
\newtheorem{corollary}[theorem]{Corollary}

\theoremstyle{definition}

\newtheorem{definition}[theorem]{Definition}

\newtheorem{remark}[theorem]{Remark}

\newtheorem*{acknowledgements}{Acknowledgments} 

\theoremstyle{remark}
\newtheorem{case}{Case}

\begin{document}
\bibliographystyle{amsalpha}
\maketitle 
\tableofcontents

\setcounter{section}{-1}
\section{Introduction}

This is the second part of our work on stable reductive varieties in
which we extend the results of the affine part
\cite{AlexeevBrion_Affine} to the projective settings.  The main aim
of this paper is to construct a compactified moduli of stable
higher--dimensional projective varieties, in a special case when
there is a nontrivial reductive group action.  Several examples of
such moduli spaces are known. The most familiar must be 
Deligne--Mumford--Knudsen moduli space $\oM_{g,n}$ of $n$--pointed genus
$g$ curves and Kontsevich's relative version for stable maps.  The
space $\oM_{g,n}$ was extended to the case of dimension two, for
surfaces $X$ of general type, and pairs $(X,D)$ of log general type by
Koll\'ar, Shepherd--Barron, the first author, and others, using the
methods of log Minimal Program
\cite{KollarShepherdBarron88,Alexeev_Boundedness,Alexeev_Mgn}.

Finally, in \cite{Alexeev_CMAV} it was extended to the case of pairs
$(X,D)$ in which the variety $X$ has arbitrary dimension but comes
with a nontrivial semiabelian group action, for example, an abelian
or toric variety. In the abelian case, this gives a moduli
compactification of $A_g$, the moduli space of principally polarized
abelian varieties, extending previous work of Namikawa and others
\cite{Namikawa_NewCompBoth}.

On the other hand, the toric case is closely related to, and
influenced by the work of Gelfand, Kapranov, Zelevinsky, Sturmfels and
others on $A$--hypergeometric functions. Each stable toric pair
$(X,D)$ comes with a canonical map to $\bP^n$ and the image is its
``shadow cycle''.  The works
\cite{KapranovSturmfelsZelevinsky_ChowPolytopes92,
  GelfandKapranovZelevinsky_Book94} describe the toric Chow variety
parameterizing these cycles. Another related moduli space is the toric
Hilbert scheme of Peeva and Stillman \cite{PeevaStillman}.

In addition to the work in the toric case, a starting point for
our investigation was M.~Kapranov's paper \cite{Kapranov} which
studies in the case of reductive group action some of the varieties
analogous to the shadow cycles.

There are several reasons why we consider the moduli problem for pairs
of varieties with divisors, rather than simply for polarized
projective varieties. The most obvious is that without the additional
structure provided by the divisor there are only finitely many
isomorphism classes of varieties of a fixed numerical type (just as
for projective toric varieties), and the varieties have infinite
automorphism groups. Hence, even though the moduli stack can still be
defined and studied, it is rather exotic, \`a la Lafforgue's ``toric
stacks''. With the divisor added, the automorphism groups become
finite. Also, just as in the case of toric or abelian varieties, a
one--parameter family may have several ``stable'' limits. With the
divisor added, the limit becomes unique, and it always exists,
possibly after a finite base change. The last, but not the least,
reason comes from the log Minimal Program. One of its lessons is that 
varieties and pairs that have good moduli spaces must have an
ample log canonical class $K_X + D$, and semi log canonical
singularities.  Since reductive varieties are rational, a
nonzero divisor is necessary.

Here is the structure of the present paper. In Section~\ref{Main
  definitions} we define polarized stable varieties and pairs. In the
next two sections we give their complete combinatorial description, in
terms of discrete data -- complexes of polytopes invariant under the
action of Weyl group, and continuous data -- certain cohomology
groups. 

In Section~\ref{Moduli of stable reductive pairs} we prove that, in
the case of multiplicity--free support, there exists a coarse moduli
space of stable reductive pairs (Theorem~\ref{thm:moduli-pairs}) and
that it is a disjoint union of projective schemes
(Theorem~\ref{thm:mod-space-projective}). Our method of proof is
different from that in the toric case  \cite{Alexeev_CMAV}.
Further, in the case of conjugation--invariant pairs we show that the
moduli space is a union of projective toric varieties corresponding to
special fiber polytopes, which we define.

In Section~\ref{sec:Connection with the log Minimal Model Program} we
prove that, as predicted by the log Minimal Program, our varieties have
semi log canonical singularities (Theorem \ref{thm:slc}).

Finally, the last Section is devoted to several generalizations.

\begin{acknowledgements}
  The first author's research was partially supported by NSF grant
  0101280.
\end{acknowledgements}

\section{Main definitions} \label{Main definitions}

We freely use notations and basic definitions from
\cite{AlexeevBrion_Affine}. In particular, $G$ will always denote a
connected reductive algebraic group over an algebraically closed
field $k$ of characteristic zero, $B,B^-$ a pair of opposite Borel
subgroups, $T$ their common torus, and $W$ the corresponding Weyl
group. Because of its fundamental importance, let us recall the
following:

\begin{definition}\label{stablereductive}
An \emph{affine stable reductive variety} (resp.~an 
\emph{affine reductive variety}) for $G$ is a connected
(resp.~irreducible) affine $\gxg$--scheme $X$ satisfying the
following conditions: 

\begin{enumerate}

\item (on singularities) $X$ is seminormal,

\item (on stabilizers) for any $x\in X$, the stabilizer
$\Stab_{\gxg}(x)$ is connected,

\item (on orbits) $X$ contains only finitely many $\gxg$--orbits, 
 
\item (group--like condition) $(B^-\times B)X^{\diag T}$ contains a
dense subset of every $\gxg$-orbit.

\end{enumerate}

\end{definition}
(Here $X^{\diag T}$ denotes the fixed point set of $\diag T\subset\gxg$.)

We now adapt this definition to the projective setting. We will
consider \emph{polarized} projective $\gxg$--varieties, that is, pairs
$(X,L)$ where $X$ is a projective scheme with $\gxg$--action, and $L$
is an ample invertible sheaf on $X$. Recall that $L$ is
\emph{linearized} if we have a lift of the action of $\gxg$ from $X$
to $L$, linear on fibers.

\begin{definition}\label{linearized}
  A \emph{linearized stable reductive variety} (resp.
  \emph{linearized reductive variety}) for $G$ consists of a connected
  (resp. irreducible) projective scheme $X$ equipped with a
  $\gxg$--action and with an ample, $\gxg$--linearized invertible
  sheaf $L$, satisfying conditions (1), (2), (3) and

\begin{enumerate}

\item[($4'$)] for any $\gxg$--orbit $\cO$, there exists a 
$T\times T$--orbit $\cO'\subseteq \cO^{\, \diag T}$
such that: $(B^-\times B)\cO'$ is dense in $\cO$, and $\diag T$
acts trivially on the fibers of $L$ at all points of $\cO'$.

\end{enumerate} 

Equivalently, $(X,L)= \left(\Proj R, \cO(1)\right)$ for a graded
algebra $R$ such that $\wX=\Spec R$ is an affine (stable) reductive
variety for the action of the group $\wG=\bG_m\times G$ (see
Theorem~\ref{thm:equiv-def}). 
\end{definition}

As an example, let $X=G/B\times G/B^-$ and let $L$ be the 
$\gxg$--linearized invertible sheaf associated with characters
$\lambda$ of $B$, and $\mu$ of $B^-$. Then $(X,L)$ is a 
polarized stable reductive variety if and only if: $\lambda$ is
regular dominant, and $\mu=-\lambda$. Other examples include
DeConcini--Procesi's wonderful compactifications of adjoint semisimple
groups; more generally, any normal projective $\gxg$--equivariant
compactification of $G$, endowed with an ample $\gxg$-linearized
invertible sheaf, is a linearized reductive variety.

\begin{definition}\label{pairs}
  A \emph{stable reductive pair} (for $G$) consists of a linearized 
  stable reductive variety $(X,L)$ and an effective ample Cartier
  divisor $D$ such that $L= \cO_X(D)$, satisfying the crucial

\begin{itemize}\label{cond-on-divisor}

\item[] (transversality condition) $\Supp D$ does not contain any
$\gxg$--orbits.

\end{itemize}
\end{definition}

\begin{remark}
  If $X$ is normal, any effective Weil divisor satisfying the
  transversality condition is automatically Cartier,
  cf. \cite{Knop94} Lemma 2.2.
\end{remark}

\begin{definition}
  We say that a stable pair $(X,D)$ is \emph{conjugation--invariant}
  if $D$ is, i.e., one has $(g,g)D=D$ for all $g\in G$.
\end{definition}

By a basic result of GIT \cite[Props. 1.4 and 1.5]{Mumford_GIT3ed} if
$X$ is normal then some positive power $L^n$ can be linearized. The
difference between assuming linearization of $L$ and of $L^n$ is a
minor technical detail which can be treated by a simple trick, see
Section~\ref{sec:Generalizations}. When studying stable pairs $(X,D)$
we can simply replace them by pairs $(X,nD)$.

A more serious situation is when no positive power of $L$ admits a
$\gxg$--linearization. This is only possible if $X$ is not normal and
the group $G$ has a nontrivial character. In this case, by
\cite[Sec.4]{Alexeev_CMAV} there exists a cover $\pi: \wX \to X$ by a
connected scheme, locally of finite type but with infinitely many
irreducible components, such that $X=\wX / \bZ^r$ and the pullback
$\pi^* L^n$ is linearizable.  This construction is absolutely crucial
for the study of degenerations of abelian varieties, but not so for
reductive groups. Hence, we also treat it as an aside in
Section~\ref{sec:Generalizations}. Until that Section, we always
assume that $L$ has a fixed linearization.

\begin{definition}
  A \emph{family of stable reductive pairs} over a $k$--scheme $S$ is
  a proper, flat morphism $\pi:X\to S$ from a scheme $X$ equipped
  with an action of the constant group scheme $\gxg$ over $S$, and with
  a relatively ample Cartier divisor $D$, such that: $L=\cO_X(D)$ is
  $\gxg$--linearized, and the fiber $(X_{\bar s},D_{\bar s})$
  over every geometric point $\bar s$ of $S$  is a stable reductive pair.

The moduli functor $\cM$ of stable reductive pairs associates to any
scheme $S$ the set $\cM(S)$ of families over $S$ modulo
isomorphisms. We will assume our schemes to be separated and
locally Noetherian $k$--schemes.
\end{definition}

\section{Polarized stable reductive varieties}
\label{Polarized stable reductive varieties}
\subsection{Classification}

Let $(X,L)$ be a polarized projective scheme, $R(X,L)$ be the ring
$\bigoplus_{n\ge0} H^0(X,L^n)$ and $\wX=\Spec R(X,L)$ be the affine
cone. Then $\wX$ is obtained by contracting the zero section of the
line bundle over $X$ corresponding to the sheaf $L\inv$, to the vertex
$0\in\wX$. Denote this contraction $\pi:\bL\inv\to \wX$, it is a
proper birational morphism. We have 
$R(X,L)= \Gamma(\bL\inv,\cO_{\bL\inv})$, and hence 
$\pi_*\cO_{\bL\inv}=\cO_{\wX}$. The multiplicative group $\bG_m$ acts
on $\wX$, and the map $\pi:\wX\setminus 0\to X$ is the geometric
quotient; it is an affine morphism with fibers the $\bG_m$--orbits,
and $\cO_X=(\pi_*\cO_{\wX\setminus 0})^{\bG_m}$. 

\begin{lemma}\label{lem:normal-seminormal}
  $X$ is normal (resp. seminormal) if and only if so is $\wX$.
\end{lemma}
\begin{proof}
  The direction from $\wX$ to $X$ follows from 
  \cite[Lem.2.1(2)]{AlexeevBrion_Affine}.
  
  If $X$ is normal then $R(X,L)$ is the ring of regular functions on
  the normal variety $\bL\inv$, and it is normal since it is the
  intersection of the integrally closed subrings $\cO_{\bL\inv,z}$ of the
  function field of $\bL\inv$.
  
  Now, assume that $X$, and hence $\bL\inv$ are seminormal and let
  $p:Z\to \wX$ be the seminormalization. The morphism 
  $\pi:\bL\inv\to \wX$ factors through $Z$ by the universal property
  of seminormalization, see f.e.
  \cite[I.7.2.3.3]{Kollar_RationalCurves96}. In the resulting exact
  sequence 
  $$ \cO_{\wX} \overset{f}{\to} p_*\cO_Z 
  \overset{g}{\to} \pi_* \cO_{\bL\inv} $$
  the composition $g\circ f$ is an isomorphism, and $g$ is
  injective. Hence, $g$ and $f$ are isomorphisms. Since $p$ is finite,
  the variety $Z$ is affine, and it follows that $p$ is an
  isomorphism.   
\end{proof}

Now, let $X$ be a $\gxg$--variety and assume that the sheaf $L$ is
$\gxg$--linearized. Then the varieties $\bL\inv$ and $\wX$ have a
natural $\wG\times\wG$--action, where 
$$
\wG = \bG_m \times G,
$$ 
with $(t_1,t_2)\in \bG_m\times\bG_m$ acting on the fibers by
multiplying by $t_1t_2\inv$. Note that 
$$
\wT = \bG_m\times T
$$
is a maximal torus of $\wG$, with character group
$$
\wLambda = \bZ\oplus\Lambda^+
$$
and positive Weyl chamber 
$$
\wLambda^+ = \bZ\times \Lambda^+.
$$

\begin{lemma}\label{lem:RXL1}
  The affine cone $\wX$ of a polarized stable reductive variety
  $(X,L)$ for the group $G$ is an affine stable reductive variety for
  the group $\wG$.
\end{lemma}
\begin{proof}
  The conditions (1) for $X$ and $\wX$ are equivalent by
  Lemma~\ref{lem:normal-seminormal}. Since the morphism
  $\wX\setminus\{0\} \to X$ is the geometric quotient by $\bG_m$, the 
  $\wG\times\wG$--orbits in $\wX$ are in bijection with the $\gxg$--orbits
  in $X$, except for the $1$--point orbit $0\in \wX$. So, the
  equivalence of conditions (2),(3) for $\wX$ and $X$ is immediate.
  Conditions (4) for $\wX$ and ($4'$) for $X$ are equivalent because
  $\diag\bG_m$ acts trivially on $\wX$.
\end{proof}

The paper \cite{AlexeevBrion_Affine} gives a complete classification
of affine stable reductive varieties for the group $\wG$ as follows:
\begin{enumerate}
\item Every affine reductive variety corresponds to an admissible cone
  in $\bR\oplus \Lambda_{\bR}$, i.e. a rational polyhedral cone
  $\sigma$ such that the relative interior $\sigma^0$ meets 
  $\bR\oplus \Lambda^+_{\bR}$ and the distinct $w\sigma^0$ ($w\in W$)
  are disjoint. The associated algebra $R_{\sigma}$ is isomorphic as
  a $\wgxg$--module to 
  $\bigoplus_{\wlambda \in \sigma \cap \wLambda^+} 
  \End V_{\wlambda}$, 
  where $V_{\wlambda}$ is a simple $\wG$--module with highest weight
  $\wlambda$. 
\item An affine stable reductive variety corresponds to a complex
  $\Sigma=\{   \sigma \}$ of admissible cones together with a
  $1$--cocycle $t$. The associated algebra $R_{\Sigma, t}$ is the
  inverse limit of the algebras $R_{\sigma}$ twisted by the cocycle $t$.
\end{enumerate}

In the case of the algebra $R=R(X,L)$, because $R_0=k$ and $R_n=0$ for
$n<0$, each $\sigma$ is the cone over a rational polytope
$(1,\delta)$, $\delta\subset \Lambda_{\bR}$.

\begin{lemma}\label{lem:RXL2}
  Each $\delta$ is a lattice polytope, i.e. its vertices are all in  
  $\Lambda$. 
\end{lemma}
\begin{proof}
  By  \cite{AlexeevBrion_Affine} Proposition 4.4, it suffices to prove
  integrality of those vertices of $\delta$ that are in
  $\Lambda^+_{\bR}$. These correspond to those orbit closures in $\wX$
  that are cones over some closed orbit $Y$ in $X$. Since $H^0(Y,L)$
  is a simple $\gxg$--module, and $L$ is globally generated (more
  generally, any nef invertible sheaf on a complete spherical variety
  is globally generated, see f.e. \cite[Thm.~4]{Brion_Inamdar94}),
  the restriction map 
  $H^0(X,L) \to H^0(Y,L)$ is surjective. The corresponding highest
  weight in $H^0(X,L)$ is of the form $(\lambda,-\lambda)$, where
  $\lambda$ is the vertex of $\delta$ associated with $Y$. So this
  vertex is integral. 
\end{proof}

Let us give the corresponding basic definitions for polytopes.

\begin{definition}\label{admissible1-proj}
A \emph{$W$--admissible polytope} is a
polytope $\delta$ in $\Lambda_{\bR}$ with vertices in $\Lambda$,
satisfying the following conditions:

\begin{enumerate}

\item The relative interior $\delta^0$ meets $\Lambda^+_{\bR}$.

\item The distinct $w\delta^0$ ($w\in W$) are disjoint. 
\end{enumerate}
We will put $\delta^+ = \delta \cap \Lambda^+_{\bR}$.
\end{definition}

Let $\sigma=\Cone \delta$ be the cone over $(1,\delta)$ in 
$\bR \oplus \Lambda_{\bR}$.  Associated with this data, we have a ring
$R_{\delta} := R_{\Cone\delta}$ and the corresponding affine reductive
variety.

\begin{definition}\label{complex-proj}
  A \emph{$W$--complex of polytopes $\Delta$ referenced by $\Lambda$}
  is a topological space $\vert\Delta\vert$ represented as a finite
  union of distinct closed subsets $\delta$ ($\delta\in\Delta$),
  together with a map $\rho:\vert\Delta\vert\to\Lambda_{\bR}$ such
  that:

\begin{enumerate}

\item $\rho$ identifies each $\delta\in\Delta$ with a polytope
in $\Lambda_{\bR}$ with vertices in $\Lambda$.

\item If $\delta\in\Delta$, then each face $\gamma\prec\delta$ is in
$\Delta$. 

\item If $\delta$, $\gamma$ in $\Delta$, then their intersection in
$\vert\Delta\vert$ is a union of faces of both.

\item $W$ acts on $\vert\Delta\vert$, the reference map $\rho$ is
$W$--equivariant, and its restriction to any subset
$\bigcup_{w\in W} w\delta$ is injective.

\end{enumerate}
(In particular, $W$ permutes the subsets $\delta$.) We will put
$|\Delta|^+ = |\Delta| \cap \Lambda^+_{\bR}.$
\end{definition}

As in the affine case, we have the set $\Delta/W$, partially ordered
by inclusion, and the complex of diagonalizable groups
$$
0\to \bigoplus_{\overline{\delta}\in\Delta/W}
\Aut^{\wgxg}(\wX_{\Cone\delta})  
\to \bigoplus_{\overline{\gamma}\prec\overline{\delta}} 
\Aut^{\wgxg}(\wX_{\Cone\gamma}) \to \cdots
$$
with the obvious differential. We denote this complex by
$C^*(\Delta/W,\Aaut)$, with cocycle groups $Z^i(\Delta/W,\Aaut)$ and
cohomology groups $H^i(\Delta/W,\Aaut)$. For each $t\in Z^1(\Delta/W,
\Aaut)$ we can define the associated algebra $R_{\Delta,t}$ and
the corresponding affine stable reductive variety, as follows.
For each $\gamma\prec\delta$, we have a natural surjective
homomorphism $R_{\delta}\to R_{\gamma}$ sending all the $\End
V_{\wlambda}$-blocks not in $\gamma$ to zero. We twist these
homomorphisms by the cocycle $t$, and define
$$ R_{\Delta,t} = \varprojlim_{\Delta/W} R_{\delta}. $$
We can now prove the following:

\begin{theorem}\label{thm:equiv-def}
  \begin{enumerate}
  \item For any polarized stable reductive variety $(X,L)$, the
    graded algebra $R(X,L)$ coincides with the algebra $R_{\Delta,t}$
    for some $W$--admissible complex of polytopes $\Delta$ and a
    1-cocycle $t$. In
    particular, all $\delta\in \Delta$ are lattice polytopes.
  \item \label{item:part2} Vice versa, for any $W$--admissible complex
    of polytopes $\Delta$ and any $1$--cocycle $t$, 
    $(X,L) := \left(R_{\Delta,t}, \cO(1) \right)$ is a polarized
    stable reductive variety. In particular, the sheaf $\cO(1)$ is
    invertible. Moreover, $R(X,L) = R_{\Delta,t}$.
  \item For any polarized stable reductive variety $(X,L)$, the sheaf
    $L$ is isomorphic to  $\cO_{\Proj R(X,L)} (1)$.
  \end{enumerate}
\end{theorem}
\begin{proof}
  (1) holds by Lemmas~\ref{lem:RXL1} and \ref{lem:RXL2}.
  
  (2) Since all vertices of polytopes $\delta$ are integral, the proof
  of Lemma \ref{lem:RXL2} shows that the zero set of
  $(R_{\Delta,t})_1$ is empty. By \cite[2.5.8]{EGA2} this implies that
  $\cO(1)$ is invertible.  The algebras $R(X,L)$ and $R_{\Delta,t}$
  coincide in degrees $n=0$ and $n\ge n_0$, they are both seminormal,
  and the corresponding affine varieties map bijectively to $\Spec
  (R_0 \oplus\oplus_{n\ge n_0} R_n)$. By uniqueness of
  seminormalization, they coincide.
  
  (3) On $\Proj R_{\Delta,t}$ one has $\cO(n)= \cO(1)^n$. Since
  $L^n=\cO(n)$ for two consecutive sufficiently large $n$, it follows
  that $L=\cO(1)$.
\end{proof}

We can now translate the classification results from the affine case
to the polarized case. 
  \cite[Prop.4.2,4.4]{AlexeevBrion_Affine} gives

\begin{theorem}\label{thm:corr-irr}
  \begin{enumerate}
  \item The polarized reductive varieties are precisely the polarized 
    varieties $(X_{\delta},L_{\delta})$, each for a unique
    $W$--admissible polytope $\delta$.
  \item The group $\Aut^{\gxg}(X_{\delta},L_{\delta})$ is
    diagonalizable with character group 
    $$\left(\wLambda\cap\lin \delta\right)/\bZ K_{\delta} =
    \bZ \oplus (\Lambda\cap\lin \delta)/\bZ K_{\delta},
    $$
    where $K_{\delta}$ denotes the set of those simple roots $\alpha$
    such that $\delta^0$ meets the hyperplane $(\alpha=0)$.
  \item The assignment $(X,L)\mapsto (X^{\diag T},L\vert_{X^{\diag T}})$
    defines a bijection from the polarized reductive varieties for
    $G$, to the polarized stable toric varieties for $T$
    with a compatible $W$--action such that the quotient by $W$ is
    irreducible.  Moreover, the $\gxg$--orbits in $X$ are in bijection
    with the $WT$--orbits in $X^{\diag T}$, and $\Aut^{\gxg}(X,L)$ is
    isomorphic to the automorphism group of the polarized
    $WT$--variety $(X^{\diag T},L\vert_{X^{\diag T}})$.
  \end{enumerate}
\end{theorem}

\cite[Prop.5.3 and Thm.5.4]{AlexeevBrion_Affine} translate to

\begin{theorem}\label{correspondence2-proj}

\begin{enumerate}
  
\item The polarized stable reductive varieties are precisely the 
  $(X_{\Delta,t},L_{\Delta,t})$, each for a uniquely determined
  $W$--admissible complex of polytopes $\Delta$ and a $1$--cohomology
  class $[t]\in H^1(\Delta/W, \Aaut)$. The irreducible components of 
  $X_{\Delta,t}$   are the varieties $X_{\delta}$ where
  $\delta\in\Delta$ is a maximal polytope 
  such that $\delta^0$ meets $\Lambda^+_{\bR}$.

\item The $\gxg$--orbits in $X_{\Delta,t}$ are in bijection with the
$W$--orbits of polytopes in $\Delta$.

\item $\Aut^{\gxg}(X_{\Delta,t},L_{\Delta,t}) = H^0(\Delta/W,\Aaut)$.
  
\item The assignment $(X,L)\mapsto (X^{\diag T}, L\vert_{X^{\diag T}})$
  defines a bijective correspondence from the polarized stable
  reductive varieties (for $G$), to the polarized stable toric
  varieties (for $T$) with a compatible action of $W$. This
  correspondence preserves orbits and automorphism groups.
\end{enumerate}
\end{theorem}

\subsection{Cohomology groups} 

\begin{theorem}\label{thm:cohomology}
Let $(X,L)$ be a polarized stable reductive variety. Then 
$H^i(X,L^n) = 0$ for all $i,n>0$.
\end{theorem}

\begin{proof}
  For irreducible $X$, the assertion follows for example
  from \cite[Cor.5.2]{Brion97}. We prove the general case by induction
  on $\dim |\Delta|$, where $(\Delta,t)$ are the combinatorial data of
  $X$.  

  Let $\delta_1,\dots,\delta_m$ be the maximal polytopes in $\Delta$
  such that $\delta_i^0\cap \Lambda^+_{\bR}\ne \emptyset$.
  Consider the following resolution
$$
0\to \cO_{X_{\Delta,t}} \to \bigoplus \cO_{\delta_{i_0}} \to 
\bigoplus \cO_{\delta_{i_0}\cap\delta_{i_1}} \to \cdots 
$$
with the differential corresponding to $t$.  Here, 
$\Delta_{i_0\dots  i_p} = \delta_{i_0} \cap \dots \cap \delta_{i_p}$ 
is a complex of polytopes. It may consist
of more than one polytope if the reference map $\rho$ is not
injective. The existence of the resolution is immediate in the affine
case, for the affine cones $\wY=\wX$, $\wX_{\delta_0\cap\delta_1}$
etc. over $Y=X$, $X_{\delta_0\cap\delta_1}$ etc. For each of these
varieties, $\pi:\wY\setminus\{0\} \to Y$ is the geometric quotient by
the $\bG_m$-action; this implies the existence of the
resolution in the projective case.

By twisting this resolution by $L^n$ and taking the hypercohomology, we
see that each $H^i(X,L^n)$ is computed by a spectral sequence with first
terms $H^q(X_{\Delta_{i_0\dots i_p}},L^n)$. Among these, all the groups
with $q>0$ are zero either by induction assumption or by the
irreducible case. Moreover, the sequence in degree $0$ splits into a
direct sum over the 
$\wlambda\in |\Cone\Delta| \cap (\id,\rho)\inv(\wLambda)$.
For each such $\wlambda$ we get an elementary 
inclusion--exclusion complex whose cohomology is $\End V_{\wlambda}$ in
degree zero, and $0$ in higher degrees. This completes the proof.
\end{proof}

\section{Pairs} \label{Pairs}

By Theorem~\ref{thm:cohomology}, each Cartier divisor $D$ with
$\cO_X(D)=L$ is given by a section
$$ 
s =\sum s_{\lambda} \in H^0(X,L) = \bigoplus \End V_{\lambda}, 
\qquad  \lambda \in |\Delta| \cap \rho\inv(\Lambda^+).
$$
The transversality condition of Definition~\ref{pairs} plainly means
that $s_{\lambda}\ne 0$ for all $\lambda$ corresponding to the closed
$\gxg$--orbits, i.e. for the vertices of all polytopes
$\delta\in\Delta$.

\begin{definition}
  The \emph{type} of a stable reductive pair $(X,D)$ is the pair
  $(\Delta,C)$, where $\Delta$ is a $W$--admissible complex of
  polytopes, and $C$ is the $W$--invariant subset of
  $|\Delta|\cap \rho\inv(\Lambda)$ defined by
  $$\lambda \in C \cap \rho\inv(\Lambda^+) \iff
  s_{\lambda}\ne 0.$$
  Then $C$ contains all vertices of polytopes $\delta\in\Delta$. We will
  put $C_{\delta}=C\cap \delta$, $C^+=C\cap \rho\inv(\Lambda^+)$, and 
  $C_{\delta}^+ = C_{\delta} \cap C^+$.
\end{definition}

It is obvious that $D$ is conjugation--invariant if and only if each
$s_{\lambda}\in\End V_\lambda$ is a scalar matrix. In this case, $s$
can be considered to be a section of the restriction of $L$ to the
stable toric variety $X^{\diag T}$, since the restriction map
$$
H^0(X,L)^{\diag G}\to H^0(X^{\diag T},L)^W
$$
is an isomorphism by 2.1 and \cite[4.5,5.4]{AlexeevBrion_Affine}. The
following theorem now readily follows:

\begin{theorem}\label{correspondence3-proj}
  The assignment 
  $\left(X,D=(s)\right) \mapsto \left( X^{\diag T},D'=(s) \right)$ 
  defines a bijective correspondence from the conjugation--invariant
  stable reductive pairs (for $G$), to the stable toric pairs (for
  $T$) with a compatible $W$--action such that the equation of $D'$
  is $W$--invariant. This correspondence preserves automorphism groups.
\end{theorem}

\begin{definition}\label{defn:Fun}
  Denote $\Fun_{\delta,C} = \Fun(C_{\delta}^+, \bZ)$, the free abelian 
  group of $\bZ$--valued functions on the finite set
  $C^+_{\delta}$. We have an exact sequence
$$
0\to \bL_{\delta,C} \to \Fun_{\delta,C} \overset{p}{\to}
(\wLambda\cap\lin\delta)/\bZ K_{\delta}  \to \bK_{\delta,C} \to 0
$$
Here the homomorphism $p$ in the middle maps $1_\lambda$ to 
$(1,\lambda) \mod \bZ K_{\delta}$, and the first and the last groups
are its kernel and cokernel respectively.
\end{definition}
Since $C^+_{\delta}$ contains all the vertices of $\delta$ in
$\Lambda^+$, the group $\bK_{\delta,C}$ is finite.

\begin{theorem}
  The set of isomorphism classes of pairs of type $(\delta,C)$ is the
  torus with character group $\bL_{\delta,C}$. The automorphism
  group of any pair of type $(\delta,C)$ is the finite diagonalizable
  group $\hbK_{\delta,C}$ with character group $\bK_{\delta,C}$.
\end{theorem}
\begin{proof}
  Indeed, the set of all pairs of type $(\delta,C)$ is a torus with
  character group $\Fun_{\delta,C}$; we just have to divide this
  set by the action of 
$$
\Aut^{\gxg}(X_{\delta},L_{\delta}) = \Hom(
  (\wLambda\cap\lin\delta)/\bZ K_{\delta} , \bG_m).
$$
\end{proof}

For a more general type $(\Delta,C)$ which corresponds to a variety
with several irreducible components, the groups $\Fun_{\delta,C}$
define a cosheaf $\Fun_{\Delta,C}$, and the tori $\Hom(
\Fun_{\delta,C}, \bG_m)$ -- a sheaf $\hFun_{\Delta,C}$.  Entirely
analogously to the toric case \cite[Sec.2]{Alexeev_CMAV}, one can
explicitly write down a complex $\hbM^*$, whose 0th (resp. 1st)
cohomology groups compute the automorphism group (resp. the set of
isomorphism classes) of stable reductive pairs of type $(\Delta,C)$.
This complex is the cone of the homomorphism $C^*(\Delta/W, \Aaut) \to
C^*(\Delta/W, \hFun_C)$. We note one important corollary:

\begin{corollary}
  The automorphism group of any stable reductive pair is finite.
\end{corollary}

Explicitly, this is the diagonalizable group $H^0(\Delta/W,
\hbK_{\Delta,C})$.

\section{Moduli of stable reductive pairs} 
\label{Moduli of stable reductive pairs}
\subsection{General remarks on families}

Let $S$ be a scheme, $H$ a reductive group, 
$\pi:X\to S$ a proper morphism with an action 
$\sigma:H_S\underset{S}{\times}X\to X$ of the constant group
scheme $H_S= H\underset{k}{\times} S$, and $F$ an $H$--linearized
coherent sheaf on $X$. For any affine open subset $S_0$ of $S$, the
$\Gamma(S_0,\cO_{S_0})$-module  
$M=\Gamma(S_0,\pi_*F)$ comes with a coaction of the group $H$.  By the
basic Lemma \cite[p.25]{Mumford_GIT3ed}, $M$ is a union of
finite--dimensional invariant subspaces. Since $\chr k=0$ and $H$ is
reductive, $M$ splits into a direct sum over the
finite--dimensional irreducible representations $V_{\lambda}$ of $H$:
$$ 
M= \bigoplus_{\lambda} M_{\lambda} \otimes_k V_{\lambda}, 
$$
for uniquely defined $\Gamma(S_0,\cO_{S_0})$-modules
$M_{\lambda}=\Hom^H(V_{\lambda},M)$. Hence the coherent sheaf $\pi_*F$
uniquely splits into a direct sum over irreducible representations.

If $H = \gxg$ and $L$ is a $\gxg$--linearized sheaf on $X$, we
consider the $\wG\times\wG$--action on it, where $\wG=\bG_m\times G$,
and $\bG_m \times \bG_m$ acts by dilation $t_1t_2\inv$ on fibers. The 
weight lattice for $\wG$ is $\wLambda=\bZ\oplus \Lambda$ and each
sheaf $\pi_*L^n$ has weight $n$ in the first variable.

\begin{theorem}\label{thm:algebras-split-as-reps}
  Let $\pi:(X,L)\to S$ be a family of polarized stable reductive
  varieties over a scheme $S$. Then the $\cO_S$--algebra
  $R(X,L)=\bigoplus_{n\ge0} \pi_*L^n$ can be written uniquely in the
  form
$$
R(X,L) = \bigoplus_{\wlambda\in \wLambda^+} 
F_{\wlambda} \otimes_k \End V_{\wlambda}
$$
where each $F_{\wlambda}$ is a locally free $\cO_S$--module of finite
rank.
\end{theorem}

\begin{proof}
  By Theorem~\ref{thm:cohomology} and Cohomology and Base Change
  (\cite{Hartshorne77} Theorem III.12.11) each sheaf $\pi_*L^n$ is a
  locally free $\cO_S$-module of finite rank.  By the remark above,
  with $H=\wgxg$ and $M=\pi_*L^n$, we have a canonical splitting
$$
M= 
\bigoplus_{\wlambda,\tilde\mu\in\wLambda^+} 
M_{\wlambda\tilde\mu}\otimes_k(V_{\tilde\mu}^*\otimes_k V_{\wlambda})
$$
Every direct summand must be a locally free sheaf of finite rank, and
by looking at any fiber we see that it is zero unless
$\wlambda=\tilde\mu$. 
\end{proof}

\begin{definition}
  The \emph{support} of a polarized stable reductive variety is
  the corresponding topological space $|\Delta|$ with its reference
  map $\rho$ to $\Lambda_{\bR}$. The variety $(X,L)$ is
  \emph{multiplicity--free} if each $\wlambda$--block 
  $\End V_{\wlambda}$ appears in $R(X,L)$ with multiplicity at most
  one. Clearly, this is equivalent to the condition that $\rho$ is
  injective. 
\end{definition}

\begin{theorem}
  Let $\pi:(X,L)\to S$ be a family of polarized stable reductive
  varieties over a connected scheme $S$. If one geometric fiber
  $(X_{\bar s},L_{\bar s})$ is multiplicity--free then any other
  geometric fiber $(X_{\bar t},L_{\bar t})$ is multiplicity--free and
  has the same support $\rho:|\Delta| \to \Lambda_{\bR}$.
\end{theorem}

\begin{proof}
  If the rank of the locally free module $F_{\wlambda}$ equals $1$
  at one $s\in S$, then the same is true at any $t\in S$.
\end{proof}

We will consider the moduli problem only for multiplicity--free
varieties. As a consequence of the last theorem, we can work with
varieties that have a fixed multiplicity--free support 
$\vert\Delta\vert\subseteq \Lambda_{\bR}$.  This will give us an open
and closed subscheme in our moduli space.

\subsection{One-parameter degenerations}
\label{sec:One-parameter degenerations}
We can work in two fairly similar situations:
\begin{enumerate}
\item $\cR_0$ is a discrete valuation ring with maximal ideal
  $(z)$, residue field $k$, and quotient field $\cK$, $S=\Spec\cR_0$
  with generic point $\eta=\Spec\cK$ and special point $s=\Spec k$;
  or 
\item $S=\Spec\cR_0$ is a smooth curve, $z\in \cR_0$ has a unique zero
  at the closed point $s\in S$, and $\cK=\cR_0[1/z]$.
\end{enumerate}
Given a family of stable reductive pairs over $S-\{s\}$, we first
extend it, after a finite base change $(S',s')\to (S,s)$, to a
family over the whole $S'$. Next, we show that such an extension is
unique. We will give the arguments for the first situation, and leave
making the obvious changes to adapt it for the second one to the
reader. 

We begin with a pair $(X_{\eta},D_{\eta})$ where
the geometric generic fiber 
$X_{\eta}\otimes_{k(\eta)} \overline{k(\eta)}$ 
is irreducible. Let $L_{\eta}=\cO_{X_{\eta}}(D_{\eta})$. After a
finite base change, $(X_{\eta},L_{\eta})$ becomes the standard polarized
reductive variety $(X_{\delta},L_{\delta})$ over $k(\eta)$ for some
$W$--admissible lattice polytope $\delta\subset \Lambda_{\bR}$. The
divisor $D_{\eta}$ is given by a section
$$
s = \sum_{\lambda\in \Lambda^+\cap\delta} s_{\lambda}, \qquad
s_{\lambda} \in \cK\otimes_k\End V_{1,\lambda}.
$$
Let $m_{\lambda} = \val_z s_{\lambda}$, so that $s_{\lambda} =
z^{m_\lambda} s'_{\lambda}$ with $s'_{\lambda}(0)\ne0$. The lower
convex envelope of the points $(w\lambda, m_{\lambda})\in
\Lambda\oplus\bZ$, $w\in \Stab_W(\delta)$, defines a $W$--admissible
height function $h:\Cone\delta\to \bR$, as in
\cite[Sec.7.4]{AlexeevBrion_Affine}.  
It is piecewise linear, $W$--invariant, and takes rational values at
all $\wlambda\in \wLambda$. After the base change $z=(z')^n$, we can
assume that $h|_{\wLambda}$ is in fact integer--valued.
Consider the subspace
$$ 
\cR= \bigoplus_{\wlambda\in \wLambda^+\cap\Cone\delta} \,
z^{h(\wlambda)}\cR_0 \otimes_k \End V_{\wlambda}
\subseteq R(X_{\eta}, L_{\eta}) = 
\bigoplus_{\wlambda\in \wLambda^+\cap\Cone\delta} \,
\cK \otimes_k \End V_{\wlambda}.
$$
By \cite[Sec.7.4]{AlexeevBrion_Affine}, $\cR$ is a subalgebra, and
$\wX = \Spec \cR$ is a family of affine stable reductive
varieties over $S$. Hence, $(X,L) = (\Proj_{\cR_0} \cR, \cO(1) )$ is a
family of polarized stable reductive varieties over $S$. Since by
definition $m_{\lambda}\ge h(1,\lambda)$, $s$ is in $\cR$, so it extends
to a regular section of $L$ and gives a compactification $D$ of
the divisor $D_{\eta}$. Moreover, since for the lower convex envelope
we have $m_{\lambda}= h(1,\lambda)$ for every vertex $\lambda$ of
$\delta$, the corresponding $\lambda$--components of $s'(0)$, the
residue of $s$ modulo $(z)$, are all nonzero. Hence, the divisor $D_0$
on the special fiber satisfies the transversality condition, and we
have constructed an extended family of stable reductive pairs over $S$.

Now consider the general case, with $X_{\eta}$ not necessarily
geometrically irreducible. After a finite base change, 
$(X_{\eta},L_{\eta})$
becomes the standard polarized variety $(X_{\Delta,t}, L_{\Delta,t})$
for some $W$--admissible complex of polytopes and $t\in Z^1(\Delta,
\Aaut)$ with coefficients in $k(\eta)$. By taking the valuations at
$z$, we can write $t=z^\gamma t'$ with $t'(0)\ne 0$. Explicitly,
$\gamma$ is a collection of homomorphisms
$\gamma_{\varphi}\in\Hom(\wLambda\cap\lin{\gamma}/\bZ K_{\gamma},\bZ)$ 
for all $\varphi\subset \delta_1\cap\delta_2$ 
which satisfy the 1--cocycle condition on triple intersections. 
Note that for any two maximal polytopes $\delta_1,\delta_2\in\Delta$
and a face $\varphi\subseteq \delta_1\cap \delta_2$, we have 
$$
s_{\delta_1}|_{X_\varphi} =
t_{\varphi}( s_{\delta_2}|_{X_{\varphi}} ).
$$
Therefore, for all faces $\varphi\subseteq\delta_1\cap \delta_2$ 
the following condition holds:
\begin{eqnarray*}
  \label{eq:h-t}
  (h_{\delta_1} - h_{\delta_2})|_\varphi = \gamma_{\varphi}.
\end{eqnarray*}
After a common base change, we can now extend each irreducible
component $(X_{\delta})_{\eta}$ of $X_{\eta}$ and glue them along
$t$. This proves the existence part of the following theorem

\begin{theorem}\label{thm:1-param-deg}
  Every stable reductive pair $(X_{\eta},D_{\eta})$ over $S-\{s\}$ can
  be extended to a family of stable reductive pairs over $(S,s)$,
  possibly after a finite base change $(S',s')\to (S,s)$. Such an
  extension is unique.
\end{theorem}

To prove the uniqueness part, we are free to make further finite
base changes. So, let $\pi:(X,D)\to S$ be a family of stable reductive
pairs over $S$. After a base change, the generic fiber becomes a
standard pair $(X, D)_{\Delta,t}$ over $\cK$. We have
$$ 
(X,L) = (\Proj_{\cR_0} R(X,L) , \cO(1) ), \qquad
R(X,L) = \bigoplus \pi_* L^n.
$$
The family $\wX=\Spec R(X,L)$ is a family of affine stable
reductive pairs over $S$. By the classification in 
\cite[7.11]{AlexeevBrion_Affine}, 
every such family corresponds to a system of integral--valued height
functions $h_{\delta}$, $\delta\in \Delta$. All we have to show is
that each of these height functions is defined by the lower envelope
of points $(\lambda, m_{\lambda}=\val_z (s_{\delta})_{\lambda} )$, 
as in our construction. But among all functions that take
integral values on $\wLambda$, this one is the only function that
satisfies the following two conditions:
\begin{enumerate}
\item $s_{\delta}$ remains regular, i.e. $h(1,\lambda) \ge m_{\lambda}$
  for $\lambda\in \delta\cap\Lambda^+$; and
\item on the special fiber, the divisor $D_0$ satisfies the
  transversality condition, i.e. $h(1,\lambda) = m_{\lambda}$ for all
  vertices of $\delta$. 
\end{enumerate}
This completes the proof of the Theorem.

\begin{definition}
  A subdivision $\Delta_h$ of a $W$--admissible polytope $\delta$
  (resp. subdivision $\Delta_{\uh}$ of a $W$--complex of polytopes)
  corresponding to a $W$--admissible height function (resp. system of
  height functions) is called \emph{coherent}.
\end{definition}
For ordinary polytopes, in absence of the Weyl group $W$, such
subdivisions are also sometimes called regular, or convex.  As a
corollary of the proof of Theorem~\ref{thm:1-param-deg}, we have

\begin{corollary}
  Combinatorially, one--parameter degenerations of stable reductive
  varieties correspond to coherent subdivisions. In the
  conjugation--invariant case, one--parameter degenerations of stable
  reductive pairs are in bijection with one--parameter degenerations
  of associated stable toric pairs with $W$--action.
\end{corollary}

\subsection{Construction of the moduli space of pairs}

Recall that the moduli functor 
$\cM: \operatorname{(Schemes)}^o \to \operatorname{(Sets)}$ 
associates to $S$ the set of isomorphism
classes of families of stable reductive pairs over $S$. We will also
consider a version of this functor with conjugation--invariant pairs,
for which the equation $s=\sum s_\lambda$ is a collection of
$\cO_S$--scalar matrices. We will denote the latter functor $\cMci$.
We are going to prove that both of these moduli functors have coarse
moduli spaces. For this, we will use the following basic tool from our
affine paper \cite{AlexeevBrion_Affine}.

Let $H$ be a connected reductive group with set of dominant weights 
$\Lambda^+_H$. Fix a function $h:\Lambda^+_H \to \bN$. Consider a flat
family $\pi:\tilde X\to S$ of affine $H$--varieties. Then we have a
decomposition 
$$
\pi_*\cO_{\tilde X} = 
\bigoplus_{\lambda\in\Lambda^+_H} F_{\lambda} \otimes_k V_{\lambda}
$$
where the $F_{\lambda}$ are flat $\cO_S$--modules. We say that $h$ is
the Hilbert function of $\tilde X$, if each $F_{\lambda}$ is locally
free of rank $h(\lambda)$. 

Fix, in addition, a finite--dimensional $H$--module $V$ and
define the moduli functor
$$
\cM_{h,V}: \operatorname{(Schemes)}^o \to \operatorname{(Sets)}
$$
by associating to $S$ the set of all closed $H$--invariant
subfamilies $\tilde X\subseteq V\times S$ with Hilbert function $h$.

\begin{theorem}[\cite{AlexeevBrion_Affine}, Thm.7.6]
  \label{affine-moduli}
  The functor $\cM_{h,V}$ is representable by a 
  scheme $M_{h,V}$ which is quasiprojective over $k$.
\end{theorem}

Here is our main existence theorem.

\begin{theorem}\label{thm:moduli-pairs}
  The functor $\cM$ (resp. $\cMci$) is coarsely represented by a
  scheme $\M$ (resp. $\Mci$) which is proper over $k$.
\end{theorem}

\begin{proof}
  Fix a type $(\Delta,C)$ and put  $\cD = |\Delta|$. There is a natural
  partial order on types: we say that $(\Delta',C')\ge (\Delta,C)$
  if $\Delta$ is a subdivision of $\Delta'$ and $C\subseteq C'$. The
  pairs of type $\ge (\Delta,C)$ clearly form an open subfunctor. If
  we prove that the latter functor is coarsely representable, then
  $\M$ will be obtained by the obvious gluing of these open
  subschemes. The idea is to represent $\cM_{\ge (\Delta,C)}$ as a
  finite quotient of an appropriate locally closed subscheme of some
  $M_{h,V}$, where the Hilbert function $h$ takes value $1$ at points
  of $\cD$, and $0$ else. We will prove the result in several steps,
  starting with the   simplest case and then adding layers of
  complexity. 
  \begin{case}
    The pairs are conjugation--invariant, and for each
    $\delta\in\Delta$ the set $C\cap\delta^+$ generates the semigroup
    $\Cone\delta^+  \cap \wLambda^+$.
  \end{case}
  Let $\pi:(X,D) \to S$ be a family of stable reductive pairs and let
  $L=\cO_X(D)$. We have
$$
R(X,L) = \bigoplus_{\wlambda\in \wLambda^+\cap \Cone\cD} 
F_{\wlambda} \otimes_{k} \End V_{\wlambda}
$$
for some invertible sheaves $F_{\wlambda}$, and this ring comes with
an element of degree one, the equation $s=\sum s_c$ of $D$. The rule 
$1_c \to s_c, \, c\in C^+$,
 defines a $\wgxg$--homomorphism of $\cO_{S}$--algebras
 \begin{eqnarray}
   \label{eq:Chow-surjection}
\cO_S\otimes_k\Sym^*( \bigoplus_{c\in C^+}\End V_{1,c}) \to R(X,L)   
 \end{eqnarray}
We claim that his homomorphism is surjective. Indeed, it is surjective
on the $\cO_S$--submodules corresponding to highest weight vectors:
when restricted to each subset $\Cone\delta^+$ this is just the
condition that the set $C\cap\delta^+$ is generating.  Since the image
is $\wgxg$--invariant and contains the highest weight vectors,
it is everything. This gives a closed $\wgxg$--invariant embedding 
$\wX\subseteq V\times S$, where $V$ is the $\wgxg$--module 
$\bigoplus_{c\in C^+}\End V_{1,c}^*$. 

We have just constructed, in a canonical way, a morphism 
$S\to M_{h,V}$ which must factor through the open subscheme $M_{h,V}^0$
parametrizing reduced schemes. Vice versa, any morphism 
$S\to M_{h,V}^0$ gives an $\cO_S$--algebra $\cR$ with an element $s$
of degree one which defines $(X,L) = (\Proj \cR, \cO(1))$ and
$D=(s)$. Hence, $M_{h,V}^0$ represents $\cMci_{\ge(\Delta,C)}$.

\begin{case}
  Conjugation--invariant pairs, arbitrary $(\Delta,C)$.
\end{case}

Starting with the sections $s_c$ of $F_{1,c}$, we would like to
construct nonzero $\cO_S$--scalar matrices 
$t_{\wlambda} \in \Gamma(S,F_{\wlambda})$ for all 
$\wlambda \in \Cone|\Delta| \cap N\wLambda^+$, 
for some fixed positive integer $N$.  We would like to do it in a
canonical way, that is, depending only on the type $(\Delta,C)$. 

Let 
$p_{\wlambda}:k.1_{\wlambda} \subseteq \End V_{\wlambda} \onto 
(\End V_{\wlambda})^{\wlambda}$ 
be the canonical isomorphism from the line of scalar matrices to the
highest weight line (of weight $(-\wlambda,\wlambda)$).  Let
$q_{\wlambda}$ be its inverse.
These morphisms extend uniquely to the $\cO_S$--module
$F_{\wlambda}\otimes_k \End V_{\wlambda}$. For any element $u= \sum
u_{\wlambda}$ of $R(X,L)$ define 
$r_{\wlambda}(u) = q_{\wlambda}( (u_{\wlambda})^{\wlambda} )$. This is
a $\cO_S$-scalar matrix in the $\wlambda$--block $\End V_{\wlambda}$
of $R(X,L)$.

Fix one polytope $\delta\in \Delta$. We are going to describe a
universal method for obtaining nowhere--vanishing scalar matrices
$t_{\wlambda}$ in our $\cO_S$--algebra for every 
$\wlambda\in N\wLambda\cap \Cone\delta^+$.

Let us start with the group algebra $k[G]$.  Suppose we have a scalar
matrix $1_{\lambda}\in \End V_{\lambda} \subset k[G]$. First of all,
we want to produce scalar matrices $t_{\mu}$ for all vertices $\mu$ of
the polytope $Q=\Conv(W\lambda) \cap \Lambda^+_{\bR}$, possibly after
replacing it with the dilated polytope $NQ$ so that these vertices
become lattice points.  We want to do it in a canonical fashion,
considering only the powers of $1_{\lambda}$ in $k[G]$ and its
components in various $\End V_{\mu}$. Now every vertex $\mu$ of $Q$ is
the centroid of a face of $\Conv(W\lambda)$ containing $\lambda$. And
by Lemma~\ref{lem:repthry} below, the $N\mu$--component of
$1_{\lambda}^N$ is nonzero for any multiple $N$ of some positive
integer $N_0$. 

Now, let $R = R_\sigma$ be the algebra of the affine reductive variety
for $G$ corresponding to a cone $\sigma$ and let $K$ be the set of
simple roots $\alpha$ such that $(\alpha=0)$ meets $\sigma^0$. Given 
$\lambda\in \sigma\cap\Lambda^+$, let 
$Q=\Conv(W_K\lambda)\cap\Lambda^+_{\bR}$ and let $\mu$ be a vertex of
$Q$. Then the previous construction still gives a nonzero scalar
matrix in the $N\mu$--block of $R$, by the description of the
multiplication in $R$ \cite[Sec.3.2]{AlexeevBrion_Affine}. This
implies that the construction still works when $R=R_{\Sigma,t}$ is the
algebra of an affine stable reductive variety if $\Sigma$ is
multiplicity-free and has a cone containing $\sigma$.

Finally, apply this construction to the graded $k$-algebra coming
from a pair of type $\ge(\Delta,C)$.
Starting with elements $s_c$ in the $(1,c)$-blocks
($c\in \operatorname{Vert}\delta$) we produce elements
$t_{(N_0,N_0 c')}$ for all $c'\in \operatorname{Vert}\delta^+$ and
some $N_0 > 0$. Once we have these elements, taking the components of
highest weight in their products gives us nonzero scalar matrices for
every $\wlambda\in N\wLambda\cap \Cone\delta^+$ for an even larger
$N$. Denote the latter matrices by $f_{\wlambda}(1_c)$ $(c\in C^+)$.

The functions $f_{\lambda}$ involve only multiplication in the algebra
$R$ and taking the components in a canonical splitting. Therefore, for
any $\cO_S$-algebra $R$ associated to a family of pairs of type
$\ge\Delta$ we can define
\begin{eqnarray} \label{eq:t}
  t_{\wlambda}= r_{\wlambda}\left(f_{\wlambda}(s_c) \right), 
  \qquad \wlambda\in N\wLambda\cap \Cone\delta^+ .
\end{eqnarray}
These are scalar $\cO_S$--matrices, and by the above they do not vanish
on any fiber.

Now, for each cone $\delta\in \Delta$ choose a finite set
$\{v=v_{i,\delta}\}$ generating the semigroup 
$\Cone\delta\cap \wLambda^+$. For each $Nv$ we have a section
$t_{Nv}$. Let us choose a scalar section $s_{v}$ in each $v$--block of
$R(X,L)$ so that  
\begin{eqnarray} 
  \label{eq:s^N=t}
r_{Nv} (s_{v}^N) = t_{Nv}
\end{eqnarray}
It is always possible to find these $N$th roots $s_v$ after a finite
base change $S'\to S$. Any two choices of $s_{v}$'s differ by $N$th
roots of unity. After the $s_{v}$'s were chosen for all
$\{v_{i,\delta}\}$, as in the first case, we have a canonical
surjection 
$$
\cO_S \otimes_k \Sym^*( \bigoplus_{v_{i,\delta}}\End V_{i,\delta}) 
\to R(X,L)
$$
which encodes our pair. The equations~(\ref{eq:t},\ref{eq:s^N=t})
define a closed subscheme $M'$ of a certain moduli scheme $M_{h,V}^0$,
and we have constructed a morphism $S'\to M'$. It is canonical up to
rescaling by $N$th roots of unity.  Clearly, the functor 
$\cMci_{\ge (\Delta,C)}$ is coarsely representable by the quotient of
$M'$ by a product of several groups of $N$th roots of unity.

\begin{case}
  General pairs, not necessarily conjugation--invariant.
\end{case}
In this case, the $s_c$ are some nonzero $n_c\times n_c$--matrices,
not scalar in general. For any of $\prod n_c^2$ choices of matrix
coefficients there is an open subfunctor in $\cM$ where the
corresponding coefficients $(s_c)_{ij}$ are invertible. For any such
choice, $(s_c)_{ij}$ defines a nonzero scalar matrix, so we can argue
as in the previous case. The $\prod n_c^2$ thus constructed schemes
glue in an obvious way (just the way a projective space is glued from
affine spaces) to a coarse moduli space $\M$.

The fact that space $M$ is proper is a consequence of
Theorem~\ref{thm:1-param-deg} and the valuative criterion of properness.
\end{proof}

\begin{lemma}\label{lem:repthry} (see also \cite{Timashev02} Lemma 1).
  Let $\lambda$ be a dominant weight and let $\mu$ be the centroid of
  a face of the polytope $\Conv(W\lambda)$ containing $\lambda$. Then
  there exists a positive integer $N_0$ such that the decomposition
  of   $V_{\lambda}^{\otimes N}$ into irreducible $G$--representations 
  contains $V_{N\mu}$, for any multiple $N$ of $N_0$.
\end{lemma}
\begin{proof}
  Consider first the case where $G$ is semisimple, and $\mu=0$. Then
  we may take $N_0 = \dim V_{\lambda}$, since the top exterior power
  $\wedge^{N_0}V_{\lambda}$ is a copy of the trivial representation
  in $V_{\lambda}^{\otimes N_0}$. 

  The case where $G$ is arbitrary and $\mu$ extends to a character of
  $G$ follows by an obvious shifting argument.

  In the general case, write 
  $\mu = \lambda -\sum_{i\in I} n_i \alpha_i$ for uniquely
  defined positive roots $\alpha_i$, and positive rational numbers
  $n_i$. Let $P\supseteq B$ be the parabolic subgroup
  associated to the subset $I$ of simple roots, and let
  $L\supseteq T$ be its Levi subgroup. Then the invariant subspace
  $V_{\lambda}^{R_u(P)}$ is a simple $L$--module with highest weight
  $\lambda$, and $\mu\in\Conv(W_L\lambda)$ is invariant under
  $W_L$. Hence it extends to a character of $L$, so that our assertion
  holds for $V_{\lambda}^{R_u(P)}$. By considering highest weight
  vectors, it also holds for $V_{\lambda}$.
\end{proof}

\subsection{Projectivity of the moduli space}
  The following fact is well known: 

  \begin{lemma}
    Let $f:M_1\to M_2$ be a quasi--finite morphism of proper
    schemes. If $M_2$ is projective then  so is $M_1$. If $f$ is
    surjective (on points) and $M_1$ is projective then so 
    is $M_2$.    
  \end{lemma}

\begin{theorem}\label{thm:mod-space-projective}
  The moduli spaces $\M$ and $\Mci$ are projective.
\end{theorem}
\begin{proof}
  For the proof, we will construct two Chow morphisms from 
  $\Mci_{\rm red}$ (resp.  $\M_{\rm red}$) to appropriate Chow
  varieties  parametrizing cycles of a certain, easily computable,
  degree in some fixed projective space $\bP$. Then the previous
  lemma, applied to the morphisms $\M_{\rm red}\to \Chow$ and
  $\M_{\rm red} \to \M$, will imply projectivity.
  
  Let $(X,D)$ be a family of conjugation--invariant pairs over $S$.
  Consider the  homomorphism of
  equation~(\ref{eq:Chow-surjection}). In general, it is not
  surjective. However, on the fiber at any $\bar s$, the image is a
  subalgebra $R'$ of $R=R(X_{\bar s},L_{\bar s})$ over which $R$ is
  finite. Indeed, the image  contains all $\wlambda$--blocks with
  $\wlambda\in N\wLambda^+\cap  |\Delta|$ for some fixed large
  positive integer $N$ depending only on $|\Delta|$. This gives us a
  finite morphism
$$
\varphi: X\to \bP\times S = 
\Proj \Sym^*( \bigoplus_{c\in C^+}\End V_{1,c}) \times S
$$
The image is not well defined as a family of subschemes of $\bP$ over
$S$. However, it gives a well--defined family of cycles over 
$S_{\rm red}$. By the universal property of the Chow variety (see f.e.
\cite[I.3.21]{Kollar_RationalCurves96}), we obtain a morphism
$f:\Mci_{\rm red}\to \Chow$. To show that $f$ is quasifinite, we have
to prove that a variety $X$ over $k$ is defined by its image in $\bP$
up to finitely many choices. $X$ is determined by the associated 
stable toric variety $Y=X^{\diag T}$. Since $\varphi$ is
$\gxg$--equivariant, $\varphi|_Y$ is the canonical morphism from $Y$
to $\bP$ given by the equation $s$ of $D$, and by
\cite[2.11.11]{Alexeev_CMAV}, the variety $Y$ is determined by
its image up to finitely many choices. This completes the proof for
$\Mci$.

For a general pair $(X,D)$ over $S$, we will construct, analogously to 
(\ref{eq:Chow-surjection}), a $\wgxg$-equivariant homomorphism from a 
bigger algebra:
$$
\Phi:\cO_S \otimes_k 
\Sym^*( \bigoplus_{c\in C^+}\End V_{1,c} \otimes_k (\End V_{1,c})^*)
\to R(X,L).   
$$
For each $\wlambda\in \Cone |\Delta^+|$, the $\wlambda$--block
$B_{\wlambda}\subset R(X,L)$ differs from 
$\cO_S \otimes_k \End V_{\wlambda}$ by an invertible sheaf.  This
implies that the sheaves $B_{\wlambda}\otimes B_{\wlambda}^*$ and 
$\cO_S \otimes_k 
\left(\End V_{\wlambda} \otimes_k (\End V_{\wlambda})^*\right)$ 
are canonically isomorphic. The element $s_c\in B_{1,c}$ defines a
linear map
$$
\cO_S \otimes_k \left( \End V_{1,c}\otimes_k (\End V_{1,c})^*\right) 
= 
B_{1,c}\otimes_{\cO_S} B_{1,c}^* \overset{(1,s_c)}{\longrightarrow}
B_{1,c}
\to R(X,L),
$$
and that gives a homomorphism $\Phi$ such that the algebra $R(X,L)$
is once again finite over its image. This homomorphism is
equivariant if we let $\gxg$ act trivially on $(\End V_{1,c})^*$.  So,
this time we get a family of cycles in a bigger projective space and a
quasifinite morphism to a Chow variety.
\end{proof}

\subsection{Structure of the moduli space}
\label{sub:Projectivity of the moduli space}

{}From Theorems~\ref{correspondence2-proj} and
\ref{correspondence3-proj}, we know that over an algebraically closed
field there exists a bijective correspondence between stable
reductive varieties, resp. pairs and $W$-equivariant stable toric
varieties, resp. pairs. This toric correspondence extends to families
over reduced schemes:

\begin{lemma}
  Let $(X,L)\to S$ be a flat family of polarized stable reductive
  varieties over a reduced scheme $S$. Then $X^{\diag T}$, considered
  with the reduced scheme structure, is flat over $S$ and 
  $(X^{\diag T},L\vert_{X^{\diag T}})$ is a family of polarized stable
  toric varieties with a compatible $W$-action.
\end{lemma}
\begin{proof}
  We can assume $S$ to be connected. By the toric correspondence over
  an algebraically closed field, on every geometric fiber we obtain a
  polarized stable toric variety with the same Hilbert polynomial. The
  flatness follows by applying a well-known criterion of flatness for
  projective morphisms over reduced Noetherian bases.
\end{proof}

As above, let $\rho:\cD\to\Lambda_{\bR}$ be a multiplicity--free
support.

\begin{corollary}
  There exists a natural bijective morphism $\Mci_{\rm red}\to
  \SP^W_{\rm red}$ from the moduli space of stable reductive pairs
  with support $\cD$ to the moduli space of stable toric pairs with
  compatible $W$--action and $W$-invariant equation, and with support
  $\cD$.
\end{corollary}

We can now apply \cite[Theorem 1.2.15]{Alexeev_CMAV} to deduce the
structure of $\Mci$. First, we need a few definitions. Given a
$W$--admissible polytope $\delta$, consider the simplex
$\alpha_{\delta,C}$ with vertices $1_c$ in the lattice
$\Fun_{\delta,C}$ of Definition~\ref{defn:Fun}. The homomorphism 
$p$ maps $\alpha_{\delta,C}$ to some polytope $q_{\delta,C}$.
\begin{definition}
  $\Sigma(\delta/W)$ will denote the fiber polytope of
  $p:\alpha_{\delta,C} \to q_{\delta,C}$ (see
  \cite{BilleraSturmfels_FiberPolytopes92} or
  \cite{GelfandKapranovZelevinsky_Book94} for the definition of fiber
  polytopes). It is a lattice polytope in $\bL_{\delta,C}$.
  Similarly, the polytope $\Sigma(\Delta/W,C)$ is defined as the
  projection of the polytope
$$
\prod_{\delta\in \Delta} \Sigma(\delta/W,C) \subset
  C_0(\Fun_{\Delta,C})
$$
to the space $H_0(\Fun_{\Delta,C}) = C_0/\partial C_1$.
\end{definition}
Now, \cite{Alexeev_CMAV} yields the following 
\begin{theorem}\label{thm:moduli-structure}
  The normalization of $\Mci_{\rm red}$ is a union of the projective
  toric   varieties corresponding to the polytopes
  $\Sigma(\Delta/W,C)$ for all  $W$--invariant subdivisions
  $(\Delta,C)$ of $(\cD,\cD\cap\Lambda)$. 
  The component containing reductive varieties is the toric variety
  for the polytope $\Sigma(\cD/W,\cD\cap\Lambda)$.
\end{theorem}

Note that the moduli space $\M$ of general pairs also has a natural
stratification by types $(\Delta,C)$ with the same inclusions as for
$\Mci$, but each stratum has a larger dimension 
$\sum_{c\in C^+} (n_c^2-1)$, where $n_c=\dim V_{1,c}$.

\section{Connection with the log Minimal Model Program}
\label{sec:Connection with the log Minimal Model Program}

\subsection{Arbitrary spherical varieties}

Let $X$ be a variety with an action of a connected reductive group
$G$. We will assume that $X$ is spherical, that is, $X$ is normal and
contains an open orbit of a Borel subgroup $B$ of $G$.  By
compactifying, we will also assume $X$ to be complete, this will not
affect our results on singularities below.  Since the open $B$--orbit
is affine, its complement has pure codimension one in $X$. Thus,
unlike in the toric case, $X$ has not one but two kinds of boundaries:
\begin{enumerate}
\item $\partial_G X$, the codimension one part of the complement
  of the open $G$--orbit in $X$,
\item $\partial_B X$, the complement of the open $B$--orbit minus
  $\partial_G X$ (which it contains).
\end{enumerate}

The following are well--known facts in the theory of spherical varieties
\cite{Brion_Inamdar94}: 

\begin{lemma} \label{lem:props-of-boundaries}
  \begin{enumerate}
  \item A canonical divisor for $X$ is $-\Delta_G - \Delta_B$, where
    $\Delta_G$ is the reduced divisor $\partial_G X$, and
    $\Delta_B$ is a unique effective divisor with support
    $\partial_B X$ (it need not be reduced).
  \item The linear system $|\Delta_B|$ has no fixed components.
  \item \label{item:spherical-resolution}
    There exists a $G$--equivariant resolution of singularities 
    $\pi:X'\to X$ such that $X'$ is a projective (spherical) variety
    satisfying:
    \begin{enumerate}
    \item $\Delta'_G$ is a divisor with normal crossings,
    \item $\Delta_B'$ is the strict preimage of $\Delta_B$,
    \item (transversality condition) \label{item:transv-Bruhat}
     $\Delta'_B$ does not contain any $G$--orbits.
    \end{enumerate}
    Then the linear system $\vert\Delta'_B\vert$ is free.
  \item \label{item:div-supp-Bruhat}
    Assuming $X$ is projective, there exists an effective ample
    divisor $H$ such that 
    $\Supp H = \partial_G X \cup \partial_B X = 
    \Supp(\Delta_G+ \Delta_B)$.
  \end{enumerate}
\end{lemma}

Recall the following 
\begin{definition}
  A pair $(X,D)$ has log canonical, resp. klt (Kawamata log terminal)
  singularities if 
  \begin{enumerate}
  \item $X$ is normal and $D$ is an effective $\bQ$--Weil divisor on
    $X$, 
  \item $K_X+ D$ is $\bQ$--Cartier,
  \item for one (and then any) log-resolution of singularities
    $\pi:X'\to X$ (i.e. $X'$ is smooth, and the exceptional locus + the
    strict preimage $\pi\inv D$ is a divisor with normal crossings), in
    the formula
    $$
    \pi^*(K_X+ D) = K_{X'} + \pi\inv D + \sum b_jE_j
    $$
    all the coefficients of irreducible divisors appearing in 
    $\pi\inv D + \sum b_j E_j$ are $\le 1$ (resp. $<1$). 
  \end{enumerate}
  The coefficients $a_j = -b_j$ are called the discrepancies of the
  pair $(X,D)$.
\end{definition}

\begin{theorem}\label{thm:arbitrary-spherical-lc}
  On any spherical variety $X$, the pair $(X,\Delta_G + |\Delta_B|)$
  has log canonical singularities, i.e. for a generic choice of a
  divisor $D_B\in |\Delta_B|$, the pair $(X, \Delta_G + D_B)$ has log
  canonical singularities.
\end{theorem}

\begin{proof}
  Consider the resolution of singularities $\pi:X'\to X$ as
  above. Then 
  $$
    \pi^*(K_X + \Delta_G + D_B) =0 =
    K_X + \Delta'_G + D_B'.
  $$
  Since $|D_B'|$ is free, for a generic choice of $D_B$, 
  $\Delta'_G +D_B'$ is a reduced divisor with normal crossings. 
\end{proof}

\begin{lemma}\label{lem:spherical-trick}
  For any spherical variety $X$, on some spherical resolution $X'$
  there exists a $\bQ$-divisor $A'$ such that $-(K_{X'}+A')$ is ample
  and the pair $(X',A')$ has klt singularities.
\end{lemma}

\begin{proof}
  For a generic $D'_B\in |\Delta'_B|$, $H'$ as in
  \ref{lem:props-of-boundaries}(\ref{item:div-supp-Bruhat}) and 
  $0< \delta \ll \varepsilon \ll 1$, the divisor
  $$
  A' = \Delta'_G + (1-\varepsilon) D'_B +
  \varepsilon \Delta'_B -\delta H'
  $$
  will do. Indeed, $K_{X'} + A' \sim -\delta H'$.  
  By the transversality condition
  \ref{lem:props-of-boundaries}(\ref{item:transv-Bruhat}) the
  divisor $\varepsilon\Delta'_B$ does not affect the discrepancies of
  the previous divisor that were equal to $-1$ (on some further log
  resolution $X''$), and by continuity in $\varepsilon$ (we are
  dealing with finitely many linear functions, after all) the new
  discrepancies are also $>-1$.  Finally, subtracting $\delta H'$
  strictly increases all $-1$ discrepancies because their locus is
  inside $\Supp \Delta'_G$.
\end{proof}

The next several statements are well known, see f.e. 
\cite{Brion_Inamdar94, Brion97}. However, the usual proofs are more
complicated, that of \cite{Brion_Inamdar94} using Frobenius splitting.
Here, we give one--line proofs using the well-known theorem, 
see f.e. \cite[Thm.1-2-5]{KMM}:

\begin{theorem}[Kawamata-Viehweg Vanishing Theorem]
  \label{KawamataViehweg}
  Let $\pi:X\to Y$ be a projective morphism from a normal variety
  $X$. Assume that the pair $(X,D)$ has klt singularities and let $L$
  be  an invertible sheaf on $X$ such that $L(-K_X-D)$ is
  $\pi$--ample. Then $R^i\pi_*L=0$ for $i>0$.
\end{theorem}

\begin{corollary}\label{cor:zero}
  For any nef line bundle $L'$ on a spherical resolution 
  $\pi:X'\to X$, one has $R^i\pi_* L' =0$ for $i\ge 1$. 
\end{corollary}

\begin{proof}
  Follows immediately from Lemma~\ref{lem:spherical-trick} and
  Theorem~\ref{KawamataViehweg}.
\end{proof}

\begin{corollary}
  Any spherical variety $X$ has rational singularities. 
\end{corollary}

\begin{proof}
  Applying the previous corollary with $L'=\cO_{X'}$ gives 
  $R^i \pi_*\cO_{X'}=0$. Moreover, $\pi_*\cO_{X'}=\cO_X$, since $X$ is
  normal. Thus, $X$ has rational singularities.
\end{proof}

\begin{corollary}
  For any nef line bundle $L$ on a complete 
  spherical variety $X$, one has $H^i(X,L) =0$ for $i\ge 1$.
\end{corollary}

\begin{proof}
  For $\pi:X'\to X$ as above, we set $L'=\pi^* L$ and apply the Leray
  spectral sequence $H^p(X,R^q \pi_* L') \follows H^{p+q}(X',L')$.  By
  Corollary~\ref{cor:zero} all $R^q \pi_* L'=0$ for $q>0$, so the
  spectral sequence degenerates in the first term and we get
  $H^p(X,L)= H^p(X',L')$,  which are zero for $p>0$ by Lemma 
  \ref{lem:spherical-trick} and Theorem \ref{KawamataViehweg}. 
\end{proof}


\subsection{Reductive varieties}

\begin{theorem}\label{thm2:reductive}
  Let $(X,D)$ be a conjugation--invariant reductive pair. Then
  for $0<\varepsilon \ll 1$ the pair $(X,\Delta_{\gxg} +
  |\Delta_{B\times B^-}| + \varepsilon D)$ has log canonical
  singularities.
\end{theorem}
\begin{proof}
  Without $\varepsilon D$, this is
  Theorem~\ref{thm:arbitrary-spherical-lc}. To prove it with
  $\varepsilon D$, we must check that $D$ does not contain a center of
  log canonical singularities of 
  $(X,\Delta_{\gxg} + |\Delta_{B\times B^-}|)$, i.e. the image of a
  divisor on a resolution $X'\to X$ that 
  has discrepancy $-1$. But any such center must be $\gxg$-invariant,
  so it must contain a $\gxg$-orbit. However, by the transversality
  condition the divisor $D$ does not contain any $\gxg$-orbits.
\end{proof}

Let $X=X_{\delta}$ be a projective reductive variety for $G$, such
that $\delta$ meets the interior of $\Lambda_{\bR}^+$. (This holds,
for example, if $X$ has maximal dimension i.e. $\dim G$.) 
The choice of opposite Borel subgroups $B,B^-$ in $G$ defines a Borel
subgroup $B\times B^-$ in $\gxg$. In this case, all coefficients in
$\Delta_{B\times B^-}$ equal $2$. (This follows from \cite{Brion97};
alternatively, this can be seen by looking at the maximal irreducible
degeneration of $X$, where the statement reduces 
to the case of $G/B$.) Hence, we have $\Delta_{B\times B^-}=2E$, where
$E=\partial_{B\times B^-} X$.  Similarly, 
$\Delta_{B^-\times B}=2E^-$.  The $T\times T$--invariant divisors
$E,E^-$ differ by a conjugation: $E^-=(w_0,w_0) E$, where $w_0\in W$
is the element exchanging $B$ and $B^-$ (that is, the longest element
of $W$). Hence, in this case we have a canonical representative
$E+E^-$ in the linear system $\vert\Delta_{B\times B^-}\vert$. 
Without the divisor $D$, we can
now prove a much stronger statement than Theorem
\ref{thm:arbitrary-spherical-lc} (it no longer holds with $\varepsilon
D$ added):

\begin{theorem}
  Let $X=X_\delta$ be a reductive variety such that $\delta$ meets the
  interior of $\Lambda^+_{\bR}$. Then the pair $(X,\Delta_{\gxg} + E +
  E^-)$ has log canonical singularities.
\end{theorem}

\begin{proof}
  It is enough to prove the statement in the most degenerate
  irreducible case: when $R(X,L)$ is a subalgebra of $R_{K,\Lambda}=
  \gr k[\wG]$ (with notation as in \cite[4.1]{AlexeevBrion_Affine}).
  Indeed, the most degenerate variety is the special fiber of the
  Vinberg family, and the original variety $X$ can be found as the fiber
  over a point in any open neighbourhood. By using the same trick as
  in the proof of Lemma~\ref{lem:spherical-trick}, $\Delta_{\gxg} + E
  + E^-$ is log canonical iff $\Delta_{\gxg} + (1-\varepsilon)(E +
  E^-) + \varepsilon \Delta_{B\times B^-} -\delta H$ is klt. But the
  property of being klt is open in families (this follows easily from
  \cite[17.6]{FAAT}, for example), and we are done.
  
  In the most degenerate case, since the stabilizer of a generic point
  is contained in $B^-\times B$, there exists a rational 
  map from $X$ to the projective variety $G/B^- \times G/B$. The
  morphism $X_1\to G/B^- \times G/B$ from the graph of
  this rational map is locally trivial, with fibers being toric
  varieties. Therefore, the singularities of $X_1$ can be resolved
  equivariantly to obtain a morphism $\pi:\wX\to G/B^- \times G/B$.
  
  Combinatorially, the first step -- degenerating -- corresponds to
  intersecting $\delta$ with the positive chamber; the second --
  taking the graph -- to cutting faces so that the new polytope is in
  the interior of $\Lambda_{\bR}^+$; and third --
  desingularization -- to cutting faces so that cones at the vertices
  become nonsingular.

  Let $f:\wX\to X$ be the resolution of singularities obtained this
  way. Since $\wX$ is also a spherical variety, one has
  $$
    K_{\wX} + \wDelta_{\gxg} + \wE + \wE^- = 0 = 
    f^*(K_X + \Delta_{\gxg} + E + E^-)
  $$
  and the divisors $\wE,\wE^-$ are the strict preimages of
  $E,E^-$. Hence, it is sufficient to prove that the singularities of
  $(\wX, \wDelta_{\gxg} + \wE + \wE^-)$ are log canonical.

  At this point we have a normally crossing divisor 
  $\wDelta_G=\partial_G \wX$. The divisor $\wE$ is the closure of 
  $$
    \text{(open $\gxg$--orbit)} \setminus
    \text{(open $B\times B^-$--orbit)},
  $$
  similarly for $\wE^-$ with the opposite orbits. Hence, both are
  pullbacks of divisors from $G/B^- \times G/B$:
  \begin{eqnarray*}
    \wE   &=& \pi^*( D^-\times G/B + G/B^- \times D), \\
    \wE^- &=& \pi^*( w_0 D^-\times G/B + G/B^- \times w_0 D),
  \end{eqnarray*}
  where $D$ is the complement in $G/B$ of the open subset 
  $B^-B/B$, and likewise for $D^-$. Thus, we are reduced to a problem
  about flag varieties.

  Let $Y= G/B$. The theorem would be proved if we could show that
  the pair $(Y,D + w_0 D)$ has log canonical singularities. For this
  purpose, we are going to use the ``Bott-Samelson resolution of
  $G/B$'', that we now recall. Let $w_0 = s_1 s_2 \cdots s_N$
  be a shortest decomposition of $w_0$ as a product of simple
  reflections; then $N= \dim G/B$. Let $P_i=B\cup Bs_iB$ be the
  corresponding minimal parabolic subgroups, then $P_i/B\simeq\bP^1$. 
  The Bott-Samelson resolution is the iterated $\bP^1$--bundle 
  $$
    Z = P_1 \overset{B}{\times} P_2 \overset{B}{\times} 
    \cdots \overset{B}{\times}  P_N/B \to G/B
  $$
  In $G/B$, the divisor $w_0D$ is the complement of 
  $w_0 B^- B/B = Bw_0B/B$. Thus, 
  $w_0D = \bigcup_i P_1 P_2 \cdots \widehat{P_i} \cdots P_N/B$, and its
  pullback to $Z$ is
  $$
    \bigcup_i P_1 \overset{B}{\times} P_2 \overset{B}{\times} \cdots
    \widehat{P_i} \cdots \overset{B}{\times} P_N/B.
    $$
    This boundary of $Z$ is a reduced divisor with normal
    crossings.  To complete the proof, we notice that the fibered
    product of two opposite Bott--Samelson resolutions 
    $$
    \wZ = Z\underset{G/B}{\times} Z^- = Z\underset{G/B}{\times} w_0 Z
    $$
    is nonsingular, and its boundary is a normal crossing
    divisor. Indeed, by Kleiman's transversality theorem
    \cite{Kleiman} there exists $g\in G$ such that  
    this property holds for the variety $Z\underset{G/B}{\times} gZ$.
    It is also clear that the set of all such $g\in G$ is open and
    invariant under left and right multiplication by $B$. Hence, it
    contains the unique open $B\times B$--coset $Bw_0B$; so it
    contains $w_0$. 
    This proves that the pair $(X,\Delta_{\gxg} + E + E^-)$ has log
    canonical singularities.

\end{proof}

\subsection{Stable reductive varieties}

Recall the following
\begin{definition}
  The pair $(X,D)$ has semi log canonical singularities if
  \begin{enumerate}
  \item $X$ is reduced and satisfies Serre's condition $S_2$,
  \item in codimension one, $X$ has simple crossings only,
  \item for the normalization $\nu:\wX\to X$, the pair 
  $(\wX,\nu\inv D + \text{double locus})$ has log canonical
    singularities. 
  \end{enumerate}
\end{definition}

Let $(X,D)$ be a stable reductive pair with support
$\rho:|\Delta|\to \Lambda_{\bR}$ such that:
\begin{enumerate}
\item $|\Delta|^+=|\Delta|\cap\rho\inv\Lambda^+_{\bR}$ is homeomorphic
  to a disc, 
\item $D$ is conjugation--invariant.
\end{enumerate}
In particular, any stable reductive variety in the connected
component of $\M^{ci}$ that contains a maximal--dimensional reductive
variety is such: in this case, $|\Delta|^+ = \delta^+$ for some
$W$--admissible polytope $\delta$.

Define $\Delta_{\gxg}$ to be the closure of union of the codimension
one orbits which correspond to the boundary of $|\Delta|$.  Hence, in
the irreducible case we get the same definition as before. For a
degeneration of an irreducible variety, we are only including the
limit of the previous boundary, and excluding the simple normal
singularities corresponding to the codimension one polytopes lying
inside $|\Delta|$. 

We may now generalize Theorem \ref{thm2:reductive} to stable reductive
varieties: 

\begin{theorem}\label{thm2:stable_reductive}
  Let $(X=X_\Delta,D)$ be a conjugation-invariant stable reductive
  pair. Then for $0<\varepsilon \ll 1$ the pair $(X,\Delta_{\gxg} +
  |\Delta_{B\times B^-}| + \varepsilon D)$ has semi log canonical
  singularities.
\end{theorem}
\begin{proof}
  By Lemma~\ref{lem:cm} below, $X$ is Cohen-Macaulay, and
  therefore $S_2$.  The normalization of $X$ is a disjoint union of
  reductive varieties $X_i$, corresponding to the polytopes in
  $|\Delta|\cap \Lambda^+_{\bR}$. The divisorial intersections of two
  irreducible components correspond to codimension one faces not lying
  in the supporting hyperplanes of $\Lambda^+_{\bR}$, and by an easy
  local computation they are generically normal crossings. By the
  transversality condition, $D$ does not contain any of these
  intersections. So, we have
$$
\text{(double locus) } +  \nu^*(\Delta_{\gxg}+
|\Delta_{B\times B^-}|
 + \varepsilon D) =
\sum (\Delta_{\gxg,i}+
|\Delta_{B_i\times B_i^-}|
+ \varepsilon D_i),
$$
and we are done by the  irreducible case.
\end{proof}

\begin{remark}
  This result confirms a prediction of the log Minimal Program that
  there should be a nice projective moduli space for pairs $(X,D)$
  such that $(X,D)$ has semi log canonical singularities and 
  $\omega_X(D)$ is ample, cf. \cite{Alexeev_Mgn}. Of course, the log
  Minimal Program remains conjectural in dimension larger than three.
\end{remark}

In the same way, we obtain:

\begin{theorem}\label{thm:slc}
  Assume that $X$ is equidimensional of maximal dimension. Then
  the pair $(X, \Delta_{\gxg} +E +E^-)$ has semi log
  canonical singularities.
\end{theorem}

Here $E, E^-$ are defined by gluing the divisors $E_i, E^-_i$
on the irreducible components.

To prove that stable reductive varieties are CM, we will prove the
stronger statement that their affine cone is CM. As for stable toric
varieties and Stanley-Reisner rings, this is a topological property: 

\begin{lemma}\label{lem:cm}
  Let $X=X_{\Sigma,t}$ be an affine stable reductive variety. Assume
  that the support $|\Sigma|'=\Sigma\cap\Lambda^+_{\bR}$ is
  homeomorphic to a convex cone. Then $X$ is Cohen-Macaulay.
\end{lemma}

\begin{proof}
  Consider the Vinberg family for $X$, as in
  \cite[Sec.7.5]{AlexeevBrion_Affine}. Since the property of being CM
  is open for fibers of this family, we can replace $X$ with its
  maximal degeneration, and assume that every cone
  $\sigma\in\Sigma$ is contained in $\Lambda^+_{\bR}$.  Then we have:
  $X^{\diag T}= W X'$ where $X'=X'_{\Sigma',t'}$ is the affine stable
  toric variety associated with a complex of cones with support
  $|\Sigma'|$. Since $\vert\Sigma'\vert$ is Cohen-Macaulay as a
  topological set, $X'$ is CM by \cite[Thm.2.3.19]{Alexeev_CMAV}.
  Moreover, $X'$ is fixed pointwise by $U\times U^-$, since we are
  dealing with the maximal degeneration. So the left and
  right stabilizers of $X'$ in $G$ are opposite parabolic subgroups
  $P$ (containing $B$) and $P^-$ (containing $B^-$). The corresponding
  set of simple roots consists of those orthogonal to $|\Sigma'|$.
  
  Let $\wX$ be the fiber product $G\times G\times^{P\times P^-} X'$.
  Then $\wX$ is CM and comes with a proper birational morphism
  $$
  \pi:\wX \to X
  $$
  (eg, $\pi$ is an isomorphism on all open $\gxg$-orbits).  By a
  result of Kempf \cite[pp.49-51]{Kempf} (consequence of the duality
  theorem for $\pi$, see e.g. \cite[Lemma 15] {BrionPolo} for
  details), $X$ is CM if the following three conditions hold:
  \begin{enumerate}
  \item $\pi_*\cO_{\wX} = \cO_X$.
  \item $R^i \pi_* \cO_{\wX}=0$ for $i \ge 1$.
  \item $R^i \pi_* \omega_{\wX}=0$ for $i \ge 1$, where $\omega_{\wX}$ 
denotes the dualizing sheaf.
  \end{enumerate}

To check these, we consider the projection 
$$
p:\wX \to G/P \times G/P^-.
$$
This is an affine morphism, and $p_*\cO_{\wX}$ is the $\gxg$-linearized
sheaf on $G/P \times G/P^-$ associated with the $P\times P^-$-module
$k[X']$. The latter decomposes as a sum of lines with weights
$(\lambda,-\lambda)$ where $\lambda\in |\Sigma'| $. Let 
$L(\lambda,-\lambda)$ be the corresponding line bundles on 
$G/P \times G/P^-$, then 
$$
H^i(\wX,\cO_{\wX}) = 
H^i(G/P \times G/P^-, \pi_*\cO_{\wX}) = 
\bigoplus_{\lambda\in|\Sigma'| } 
H^i(G/P \times G/P^-, L(\lambda,-\lambda)).
$$ 
Using the Borel--Weil--Bott theorem, it follows that 
$H^0(\wX,\cO_{\wX}) = k[X]$ and $H^i(\wX,\cO_{\wX}) = 0$ for $i\ge 1$.
Since $X$ is affine, this implies (1) and (2).

For (3), note that
$$
\omega_{\wX} = p^* \omega_{G/P\times G/P^-} \otimes \omega_p,
$$
where $\omega_p$ denotes the relative dualizing sheaf. Therefore,
$$
p_*\omega_{\wX} = \omega_{G/P\times G/P^-} \otimes p_*\omega_p.
$$
This is an equality of $\gxg$-linearized sheaves, and $p_*\omega_p$ 
is the $\gxg$-linearized sheaf on $G/P \times G/P^-$ associated with
the $P\times P^-$-module \newline $H^0(X',\omega_{X'})$.

We claim that this module decomposes as a sum of lines with weights
$(\lambda,-\lambda)$ where $\lambda\in |\Sigma'| ^0$ (the relative
interior of $|\Sigma'|$). Such $\lambda$'s correspond to ample line
bundles on $G/P \times G/P^-$, so that
$H^i(G/P \times G/P^-, p_*\omega_{\wX}) = 0$ for $i\ge 1$, by the
Kodaira vanishing theorem. Hence the claim completes the proof.

To prove the claim, we have to show that $H^0(X',\omega_{X'})$ is the
ideal $I_{X'}$ of $\Delta_{\gxg}$ in $k[X']$. Indeed, the
codimension--$1$ singularities of $X'$ correspond to codimension--$1$
cones $\sigma'\in\Sigma'$ that intersect $|\Sigma'|^0$. By an  
easy local computation they are generically normal crossings. Hence,
if $Z$ is the closed subset corresponding to the union of
codimension--$2$ cones then $X'_0=X'\setminus Z$ is Gorenstein.  Let
$i:X'_0\to X'$ be the inclusion. We have $\omega_{X'_0} = I_{X'_0}$,
as a combination of two results: for toric varieties, and for the
curve $\Spec k[x,y]/(xy)$, with $\deg x=1$, $\deg y=-1$.

Since the sheaf $\omega_{X'}$ is torsion free, there is an exact
sequence 
$$
0 \to \omega_{X'} \to i_*\omega_{X'_0}=i_*I_{X'_0}= I_{X'}
\to Q \to 0
$$
and the quotient $Q$ has support in codimension $\ge2$. Since $X'$ is
Cohen-Macaulay, the sheaf $\omega_{X'}$ is Cohen-Macaulay, and in
particular $S_2$. Hence, $\Ext^1(Q,\omega_{X'})=0$ and the sequence
splits. But the sheaf in the middle is torsion free, so $Q=0$. 

\end{proof}

\section{Generalizations}\label{sec:Generalizations}
\subsection{Polarized varieties and pairs without linearization}

We call a polarized $\gxg$-variety $(X,L)$ \emph{potentially
  $\gxg$--linearizable} if some positive power of $L$ is linearizable.
Equivalently, $L$ is $G_1\times G_1$--linearizable for a finite
central extension $G_1$ of $G$.
The affine cone $\wX$ of such a variety has an action of a group $\wG$
which is a central $\bG_m$--extension $\wG$ of $G$ (as opposed to the
trivial extension $\bG_m\times G$). That is the only change we need to
make, after which our classification results and construction of the
moduli go through.

For a non--normal but seminormal $X$ the sheaf $L$ is $G_1\times
G_1$--linearizable if the action of the connected center is
linearizable. That may fail, the simplest example being a plane nodal
cubic curve with the action of $\bC^*$ and $L=\cO(1)$. However,
by the result of \cite[Sec.4]{Alexeev_CMAV} there exists an infinite
$\bZ^r$--cover $\wX$ of $X$ so that the pull--back of $L$ is
linearizable. This cover is connected but is only locally of finite
type (in the simplest example above, it is an infinite chain of
$\bP^1$'s). In this case, we can classify the varieties and pairs, and
construct their moduli, by classifying the infinite covers first, and
then adding the $\bZ^r$-action. Formally, the answers are the same but
the $W$--complexes get replaced by $W\times\bZ^r$--complexes of
polytopes. We did not have a compelling reason to write down the gory
details, since the bookkeeping gets very tedious. But this could be
easily done, following the toric case of \cite{Alexeev_CMAV}.

\subsection{Pairs $(X,D_1,\dots,D_n)$ with several divisors}

In complete analogy with $n$-pointed genus $g$ stable curves and their
moduli space $\M_{g,n}$, one can consider pairs with $n$ divisors:

\begin{definition}
  A stable pair with $n$ divisors is a polarized stable reductive
  variety $(X,L)$ together with $n$ Cartier divisors $D_1\dots D_n$
  such that $L=\cO_X(\sum D_i)$ and the pair $(X, \sum D_i)$ is a
  stable pair in our previous definition.
\end{definition}

\begin{theorem}  For each $n$, there exists a coarse moduli space of
  stable pairs with $n$ divisors, having multiplicity-free
  support. Every connected component of this space is projective. 
\end{theorem}
\begin{proof}
  From the transversality condition on $D$ and the group action, it
  follows that each $D_i$ is nef. Hence, for a fixed numerical type of
  $D$ there are only finitely many possibilities for the $D_i$'s. The
  moduli space is now constructed as a closed subscheme in the moduli
  space of pairs $(X,D)$ consisting of those pairs that split into $n$
  parts. 
\end{proof}

\begin{remark}
One can give a combinatorial description of this space, similar to
Theorem~\ref{thm:moduli-structure}. Instead of the fiber polytope
$\Sigma( \alpha_{\delta} \to q_{\delta} )$, the component containing
reductive varieties will be the toric variety for the polytope
$\Sigma(\prod_{i=1}^n \alpha_{\delta_i} \to q_{\delta})$, etc.
\end{remark}




\end{document}